\documentclass[11pt,a4paper]{article}
\usepackage{indentfirst}
\usepackage{amsmath,amssymb,amsfonts,amsthm,graphics}
\usepackage{makeidx}
\usepackage{color}
\usepackage[dvips]{graphicx}
\usepackage{eufrak}
\usepackage{mathrsfs}

\begin{document}
\title{\Large \bf Additive codes over $GF(4)$ from circulant graphs\footnote{Supported by ``973" program No.2013CB834204.}}
\author{\small Xueliang Li, Yaping Mao, Meiqin Wei\\
\small Nankai University, Tianjin 300071, China\\
\small Center for Combinatorics and LPMC-TJKLC\\
\small E-mails: lxl@nankai.edu.cn; maoyaping@ymail.com; weimeiqin@mail.nankai.edu.cn\\
\small and\\
\small Ruihu Li\\
\small The air force engineering University\\
\small Institute of science, Xi'an 710051, China\\
\small E-mail: liruihu@aliyun.com}
\date{}
\maketitle
\begin{abstract}
In $2006$, Danielsen and Parker \cite{DP} proved that every
self-dual additive code over $GF(4)$ is equivalent to a graph code.
So, graph is an important tool for searching (proposed) optimum
codes. In this paper, we introduce a new method of searching
(proposed) optimum additive codes from circulant graphs.\\[2mm]


\end{abstract}

\section{Introduction}

We define $GF(4)=\{0,1,\omega,\omega^{2}\}$, where
$\omega^{2}=1+\omega$. An additive code $C$ over $GF(4)$ of length
$n$ is an additive subgroup of $GF(4)$. Clearly, $\mathcal {C}$
contains codewords for some $0\leq k\leq 2n$, and can be defined by
a $k\times n$ generator matrix, with entries from $GF(4)$, whose
rows span $\mathcal {C}$ additively. We call $\mathcal {C}$ an
$(n,2^{k})$ code. The \emph{Hamming weight} of $u\in GF(4)$, denoted
by $wt(u)$, is the number of non zero components of $u$. The
\emph{Hamming distance} between two vectors $u=(u_1,\cdots,u_n)$,
$v=(v_1,\cdots,v_n)$ is $wt(u-v)$. The \emph{minimum distance} $d$
of a code is defined as the smallest possible distance between pairs
of distinct codewords. The \emph{conjugation} of $x\in GF(4)$ is
defined by $\bar{x}=x^{2}$, and the \emph{trace map} is defined by
$Tr(x)=x+\bar{x}$. The \emph{Hermitian trace inner product} of
$u=(u_1,\cdots,u_n)$ and $v=(v_1,\cdots,v_n)$, with
$u_i,v_j\in{GF(4)}$, is given by $u\ast v=Tr(u\cdot
\bar{v})=\sum\limits_{i=1}^{n}Tr(u_{i}\bar{v_{i}})=\sum\limits_{i=1}^{n}(u_{i}v_{i}^{2}+u_{i}^{2}v_{i})$.
We define \emph{the dual of the code} $\mathcal {C}$ with respect to
the Hermitian trace inner product, $\mathcal {C}^{\bot}=\{u\in
GF(4)\ |\ u\ast c=0\ for\ all\ c\in \mathcal {C}\}$. Then $\mathcal
{C}$ is \emph{self-orthogonal} if $\mathcal {C}\subseteq \mathcal
{C}^{\bot}$, and $\mathcal {C}$ is \emph{self-dual} if $\mathcal
{C}=\mathcal {C}^{\bot}$. A graph code is an additive code over
$GF(4)$ that has a generator matrix of the form $\mathcal
{C}=\Gamma+\omega I$, where $I$ is the identity matrix and $\Gamma$
is the adjacency matrix of a simple undirected graph.

A code is called \emph{optimum} if it meets both lower and upper
bounds in the Code Tables, and a \emph{proposed optimum} code if it
only meets the lower bound in the Code Tables. The distribution of a
code is the sequence $(A_{0},A_{1},\cdots,A_{n})$, where $A_{i}$ is
the number of codewords of weight $i$. The \emph{weight enumerator}
of the code is the polynomial
$$
W(z)=\sum\limits_{i=0}^{n}A_{i}z^{i}.
$$

Let us now introduce some concepts and notions from Graph Theory. An
\emph{undirected graph} $\Gamma=(V,E)$ is a set
$V(\Gamma)=\{v_1,v_2,\cdots,v_n\}$ of vertices together with a
collection $E(\Gamma)$ of edges, where each edge is an unordered
pair of vertices. The vertices $v_i$ and $v_j$ are \emph{adjacent}
if $\{v_i,v_j\}$ is an edge. Then $v_j$ is a \emph{neighbour} of
$v_{i}$. All the neighbours of vertex $v_{i}$ in graph $\Gamma$ form
the \emph{neighbourhood} of $v_{i}$, and it is denoted by
$N_{\Gamma}(v_{i})$. The \emph{degree of a vertex} $v$ is the number
of vertices adjacent to $v$. A graph is \emph{k-regular} if all
vertices have the same degree $k$. The \emph{adjacency matrix}
$A=(a_{ij})$ of a graph $\Gamma=(V,E)$ is a symmetric $(0,1)$-matrix
defined as follows: $a_{i,j}=1$ if the $i$-th and $j$-th vertices
are adjacent, and $a_{i,j} =0$ otherwise.

\begin{figure}[h,t,b,p]
\begin{center}
\scalebox{0.6}[0.6]{\includegraphics{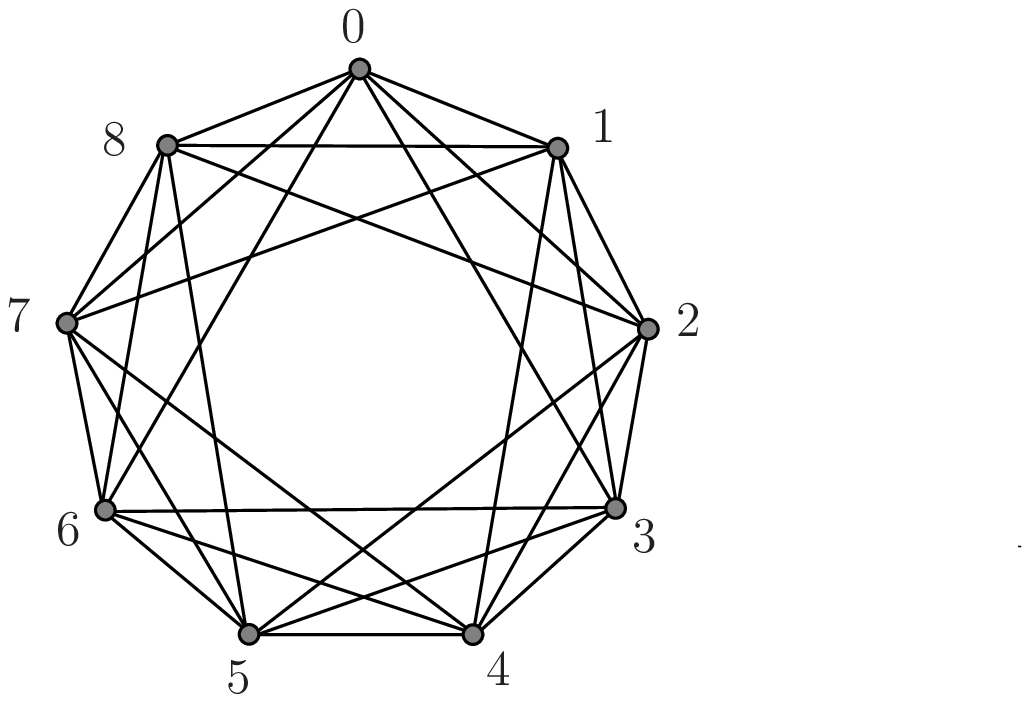}}\\
Figure 1: The circulant graph $C(9,\{1,2,3\})$ \label{fig1}
\end{center}
\end{figure}

Circulant graphs and their various applications are the objects of
intensive study in computer science and discrete mathematics, see
\cite{Bermond, Boesch, Mans, Muzychuk}. Recently, Monakhova
published a survey paper on this subject, see \cite{Monakhova}. Let
$S=\{a_{1},a_{2},\cdots,a_{k}\}$ be a set of integers such that
$0<a_{1}<\cdots<a_{k}<\frac{n+1}{2}$ and let the vertices of an
$n$-vertex graph be labelled $0,1,2,\cdots,n-1$. Then the
\emph{ciculant graph} $C(n,S)$ has $i\pm a_{1},i\pm
a_{2},\cdots,i\pm a_{k} \ (mod\ n)$ adjacent to each vertex $i$. A
\emph{circulant matrix} is obtained by taking an arbitrary first
row, and shifting it cyclically one position to the right in order
to obtain successive rows. We say that a circulant matrix is
generated by its first row. Formally, if the first row of an
$n$-by-$n$ circualant matrix is $a_{0},a_{1},\cdots,a_{n-1}$, then
the $(i,j)^{th}$ element is $a_{j-i}$, where subscripts are taken
modulo $n$. The term circulant graph arises from the fact that the
adjacency matrix for such a graph is a circulant matrix. For
example, Figure \ref{fig1} shows the circulant graph
$C(9,\{1,2,3\})$.

In $2002$, Tonchev \cite{Tonchev} first set up a relationship
between a linear binary code and the adjacency matrix of an
undirected graph. In $2006$, Danielsen and Parker \cite{DP} proved
that every self-dual additive code over $GF(4)$ is equivalent to a
graph code. Recently, finding optimum codes from graphs has received
a wide attention of many researchers, see \cite{Danielsen,
Danielsen2, Danielsen3, Danielsen4, DP, GGMG, HP, Tonchev,
Varbanov}. In \cite{DP}, Danielsen and Parker showed that two codes
are equivalent if and only if the corresponding graphs are
equivalent with respect to local complementation and graph
isomorphism. They used these facts to classify all codes of length
up to $12$. In 2012, Danielsen \cite{Danielsen} focused his
attention on additive codes over $GF(9)$ and transformed the problem
of code equivalence into a problem of graph isomorphism. By an
extension technique, they classified all optimal codes of lengths
$11$ and $12$. In fact, a computer search reveals that circulant
graph codes usually contain many strong codes, and some of these
codes have highly regular graph representations, see
\cite{Varbanov}. In \cite{Danielsen}, Danielsen obtained some
optimum additive codes from circulant graphs in 2005. Later,
Varbanov investigated additive circulant graph codes over $GF(4)$,
see \cite{Varbanov}.

In this paper, we find out some optimum additive codes from some
special circulant graphs. This paper is organized as follows. In
Section $2$, we propose a new method to find additive optimum codes
from circulant graphs. Inspired by the optimum additive codes
obtained by Danielsen \cite{Danielsen}, we focus on the dense
circulant graphs and get some optimum additive codes in Section $3$.

\section{New codes from sparse circulant graphs}

In fact, Glynn et al. \cite{GGMG} obtained an optimum code from a
circulant graph, called 5-valent graph. Recall that $\Gamma_{12}$ is
a circulant graph of order $12$; see Figure 1 $(a)$.

\begin{figure}[h,t,b,p]
\begin{center}
\scalebox{0.6}[0.6]{\includegraphics{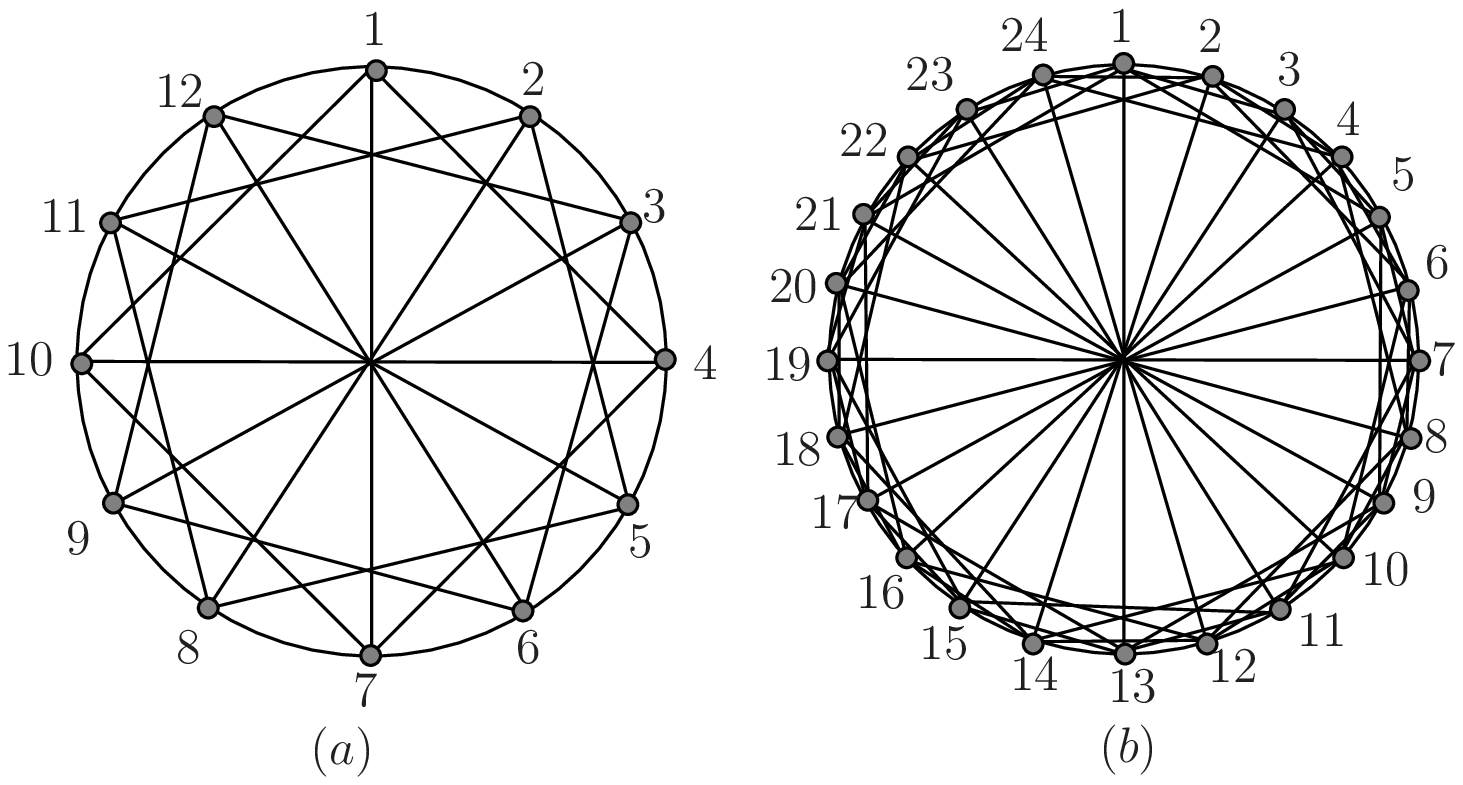}}\\
Figure 1: $(a)$ a $5$-valent graph $\Gamma_{12}$; $(b)$ the
circulant graph $\Gamma_{24}$.
\end{center}
\end{figure}

Let $V(\Gamma_{12})=\{u_1,u_2,\cdots,u_{12}\}$. For the vertex
$u_1$, we let $E_1=\{v_1v_2, v_1v_4, v_1v_7, v_1v_{10},
v_1v_{12}\}\subseteq E(\Gamma_{12})$. For the vertex $u_2$, we just
rotate the above vertices and edges, that is, we only permit the
existence of the set
$E_2=\{v_2v_3,v_2v_4,v_2v_6,v_2v_{10},v_2u_{12}\}\subseteq E(G)$ of
edges. For each vertex $u_i\in V(\Gamma)\setminus
\{u_1,u_2\}=\{u_3,u_4,\cdots,u_{12}\}$, we can also obtained the
sets $E_i \ (3\leq i\leq 17)$ of edges. Observe that
$E(\Gamma_{12})=\bigcup_{i=1}^{12}E_i$.

The adjacency matrix of the graph $\Gamma_{12}$ is the following
circulant matrix.
\begin{equation*}
A_{12}=\left(
\begin{array}{ccccccccccccccccc}
0&1&0&1&0&0&1&0&0&1&0&1\\
1&0&1&0&1&0&0&1&0&0&1&0\\
0&1&0&1&0&1&0&0&1&0&0&1\\
1&0&1&0&1&0&1&0&0&1&0&0\\
0&1&0&1&0&1&0&1&0&0&1&0\\
0&0&1&0&1&0&1&0&1&0&0&1\\
1&0&0&1&0&1&0&1&0&1&0&0\\
0&1&0&0&1&0&1&0&1&0&1&0\\
0&0&1&0&0&1&0&1&0&1&0&1\\
1&0&0&1&0&0&1&0&1&0&1&0\\
0&1&0&0&1&0&0&1&0&1&0&1\\
1&0&1&0&0&1&0&0&1&0&1&0\\
\end{array}\right)
\end{equation*}
The above matrix can also be obtained by the following vector
$$
\alpha_{12}=(0,1,0,1,0,0,1,0,0,1,0,1).
$$
Observe that this vector just corresponds to the set $E_1$ of edges,
which is an expression of the adjacent relation about the vertex
$u_1$. We conclude that a $5$-valent graph can be determined by the
edge set $E_1$, and the adjacency matrix of this graph is determined
by the above vector $\alpha_{12}$.

Furthermore, the matrix $A_{12}+\omega I$ is
\begin{equation*}
A_{12}'=A_{12}+\omega I=\left(
\begin{array}{ccccccccccccccccc}
\omega&1&0&1&0&0&1&0&0&1&0&1\\
1&\omega&1&0&1&0&0&1&0&0&1&0\\
0&1&\omega&1&0&1&0&0&1&0&0&1\\
1&0&1&\omega&1&0&1&0&0&1&0&0\\
0&1&0&1&\omega&1&0&1&0&0&1&0\\
0&0&1&0&1&\omega&1&0&1&0&0&1\\
1&0&0&1&0&1&\omega&1&0&1&0&0\\
0&1&0&0&1&0&1&\omega&1&0&1&0\\
0&0&1&0&0&1&0&1&\omega&1&0&1\\
1&0&0&1&0&0&1&0&1&\omega&1&0\\
0&1&0&0&1&0&0&1&0&1&\omega&1\\
1&0&1&0&0&1&0&0&1&0&1&\omega\\
\end{array}\right)
\end{equation*}
This matrix is also a circulant matrix, which can be obtained by the
following vector
$$
\alpha_{12}'=(\omega,1,0,1,0,0,1,0,0,1,0,1).
$$
From the matrix $A_{12}+\omega I$, we can get a graph code $\mathcal
{C}_{12}$. From the Code Tables, we know that $(12, 2^{12},6)$ is an
optimum additive code over $GF(4)$. Then $n=12$ and $d_{12}=6$,
where $d_{12}$ is the minimum distance of the code $\mathcal
{C}_{12}$.

The above statement suggests the following method for finding
optimum codes.

\textbf{Step 1.} Given an even integer $n$. Denote by $L_n$ the
lower bound of the additive code ($n$,$2^{n}$) over $GF(4)$. From
Code Tables, we find the exact value of $L_n$ for the given $n$. We
now construct a circulant graph $\Gamma_n$ by a set $E_1$ of edges
as follows.

\textbf{Step 1.1.} Arrange the vertices from
$V(\Gamma_n)=\{u_1,u_2,\cdots,u_n\}$ in a circular order.

\textbf{Step 1.2.} Determine the set $E_1$ of edges satisfying
$|E_1|=L_n-1$ or $|E_1|=L_n+1$, where $E_1=N_{\Gamma_n}(u_1)$. If
$|E_1|=L_n+1$, then
\begin{eqnarray*}
E_1&=&\{u_1u_{2},u_1u_{3}\}\cup \{u_{1}u_{3+2\cdot
1},u_{1}u_{3+2\cdot 2},\cdots, u_{1}u_{3+2\cdot
\frac{L_n-4}{2}},\}\cup
\{u_1u_{\frac{n}{2}+1}\}\\
& &\cup \{u_{1}u_{n-1-2\cdot
\frac{L_n-4}{2}},\cdots,u_{1}u_{n-1-2\cdot 2},u_{1}u_{n-1-2\cdot
1}\}\cup \{u_1u_{n},u_1u_{n-1}\}\\
&=&\{u_1u_{2},u_1u_{3}\}\cup \{u_{1}u_{5},u_{1}u_{7},\cdots,
u_{1}u_{L_n-1}\}\cup
\{u_1u_{\frac{n}{2}+1}\}\\
& &\cup \{u_{1}u_{n-L_n+3},\cdots,u_{1}u_{n-5},u_{1}u_{n-3}\}\cup
\{u_1u_{n},u_1u_{n-1}\}
\end{eqnarray*}
If $|E_1|=L_n-1$, then
\begin{eqnarray*}
E_1&=&\{u_1u_{2},u_1u_{3}\}\cup \{u_{1}u_{3+2\cdot
1},u_{1}u_{3+2\cdot 2},\cdots, u_{1}u_{3+2\cdot
\frac{L_n-6}{2}},\}\cup
\{u_1u_{\frac{n}{2}+1}\}\\
& &\cup \{u_{1}u_{n-1-2\cdot
\frac{L_n-6}{2}},\cdots,u_{1}u_{n-1-2\cdot 2},u_{1}u_{n-1-2\cdot
1}\}\cup \{u_1u_{n},u_1u_{n-1}\}\\
&=&\{u_1u_{2},u_1u_{3}\}\cup \{u_{1}u_{5},u_{1}u_{7},\cdots,
u_{1}u_{L_n-3}\}\cup
\{u_1u_{\frac{n}{2}+1}\}\\
& &\cup \{u_{1}u_{n-L_n+5},\cdots,u_{1}u_{n-5},u_{1}u_{n-3}\}\cup
\{u_1u_{n},u_1u_{n-1}\}
\end{eqnarray*}

\textbf{Step 2.} By the edge set $E_1$, we write the vector
$\alpha_n$ corresponding to $E_1$. If $|E_1|=L_n+1$, then

\quad \quad \quad \ \ $u_{2} \ u_{3}$ \quad \ $u_{5}$\quad \quad
\quad \ $u_{L_n-1}$ \quad \quad \quad  \quad \quad \
$u_{\frac{n}{2}+1}$ \quad \ \quad \quad \quad $u_{n-L_n+3}$ \quad \
\ $u_{n-3}$\quad $u_{n-1} \ u_{n}$
\begin{equation*}\alpha_n=
\left(
\begin{array}{lllllllllllllllllllllllll}
0&1&1&0&1&\cdots&0&1&0&0&\cdots&0&1&0&0&\cdots&0&1&0&\cdots&1&0&1&1\\
\end{array}\right)
\end{equation*}
If $|E_1|=L_n-1$, then

\quad \quad \quad \ \ $u_{2} \ u_{3}$ \quad \ $u_{5}$\quad \quad
\quad \ $u_{L_n-3}$ \quad \quad \quad  \quad \quad
$u_{\frac{n}{2}+1}$ \quad \quad \quad \quad \quad $u_{n-L_n+5}$
\quad \ $u_{n-3}$\quad $u_{n-1} \ u_{n}$
\begin{equation*}\alpha_n=
\left(
\begin{array}{lllllllllllllllllllllllll}
0&1&1&0&1&\cdots&0&1&0&0&\cdots&0&1&0&0&\cdots&0&1&0&\cdots&1&0&1&1\\
\end{array}\right)
\end{equation*}

\textbf{Step 3.} Change the first component of the vector $\alpha_n$
to $\omega$. Denote by $\alpha'_n$ the new vector. We generate a
circulant matrix $A'_n$ from $\alpha'_n$.
$$
\alpha_n'=(\omega,1,1,0,1,0,1,\cdots,0,1,0,0,\cdots,0,1,0,0,\cdots,0,1,0,1,0,\cdots,1,0,1,1)
$$

\textbf{Step 4.} By Algorithm $1$, we obtain the minimum distance
$d_n$ of the code $\mathcal {C}_n$ and determine whether
$d_{n}=L_n$. If so, the code $\mathcal {C}_n$ is a proposed optimum
code.

Below is an algorithm (running in SAGE). For more details, we refer
to \cite{Stein}.

\begin{tabbing}
\noindent\rule[0.25\baselineskip]{\textwidth}{2pt}
\\\textbf{Algorithm} 1: Minimum distance of a circulant graph
code\\
\noindent\rule[0.25\baselineskip]{\textwidth}{1pt}\\
Input: the value of $n$, the generator vector $\alpha_n$ of a
circulant graph code $\mathcal{C}_n$\\
Objective: the minimum distance of the circulant graph
code $\mathcal{C}_n$\\
1. input the value of $n$, the generator vector $\alpha_n=(b_1,b_2,\cdots,b_n)$;\\

2. obtain the generator matrix $G$ of the circulant graph
code $\mathcal{C}_n$; \\

3. get the minimum distance of the circulant graph
code $\mathcal{C}_n$.\\[0.2cm]

Take an example, let $n=24$ and
$\alpha_n=(\omega,1,1,0,1,0,0,0,0,0,0,0,1,0,0,0,0,0,0,0,1,0,1,1)$.\\
The algorithm details are stated as follows:\\[0.1cm]

\emph{Program}:\\
~~~~~~$F.<x>=GF(4,'x')$\\
~~~~~~$n=24$;\\
~~~~~~$a=[x,1,1,0,1,0,0,0,0,0,0,0,1,0,0,0,0,0,0,0,1,0,1,1]$\\
~~~~~~$m=matrix(F,n,n,[[a[(i-k)\%n]$ for $i$ in $[0..(n-1)]]$ for $k$ in $[0..n-1]])$;\\
~~~~~~$f=lambda \ s:sum(map(lambda \ x:m[x],s))$;\\
~~~~~~$s=[ \ ]$;\\
~~~~~~for $k$ in $[1..8]$:\\
~~~~~~~~~~~$t=min([n-list(i).count(0)$ for $i$ in Subsets $(range(n),k).map(f)])$;\\
~~~~~~~~~~~$s+=[t]$;\\
~~~~~~~~~~~print $s$;\\[0.2cm]

\emph{Output}:~[8]\\

~~~~~~~~~~~[8,\ 8]\\

~~~~~~~~~~~[8,\ 8,\ 10]\\

~~~~~~~~~~~[8,\ 8,\ 10,\ 8]\\

~~~~~~~~~~~[8,\ 8,\ 10,\ 8,\ 8]\\

~~~~~~~~~~~[8,\ 8,\ 10,\ 8,\ 8,\ 8]\\

~~~~~~~~~~~[8,\ 8,\ 10,\ 8,\ 8,\ 8,\ 8]\\

~~~~~~~~~~~[8,\ 8,\ 10,\ 8,\ 8,\ 8,\ 8,\ 10]\\

\emph{Result}:~~The minimum element of the last array is the minimum
distance $d_{24}$ of
the code $\mathcal{C}_{24}$,\\
 that is, $d_{24}=8$. \\

\noindent\rule[0.25\baselineskip]{\textwidth}{1pt}
\end{tabbing}

Inspired by the graph code $\mathcal {C}_1$ which corresponds to the
$5$-valent graph, we hope to find out some other optimum additive
codes for $n=24$.

\textbf{Step 1.} Recall $L_{24}$ is the lower bound of the additive
code $(24,2^{24})$ over $GF(4)$. From Code Tables, we find the exact
value of $L_{24}$, $L_{24}=8$. We now construct a circulant graph
$\Gamma_{24}$ by a set $E_1$ of edges as follows.

\textbf{Step 1.1.} Arrange the vertices from
$V(\Gamma_n)=\{u_1,u_2,\cdots,u_{24}\}$ in a circular order.

\textbf{Step 1.2.} Determine the set $E_1$ of edges satisfying
$|E_1|=L_{24}-1=7$ or $|E_1|=L_{24}+1=9$, where
$E_1=N_{\Gamma_n}(u_1)$. If $|E_1|=9$, then
$$
E_1=\{u_1u_{2},u_1u_{3},u_{1}u_{5},u_{1}u_{7},u_1u_{13},u_{1}u_{19},u_{1}u_{21},u_1u_{23},u_1u_{24}\}
$$

If $|E_1|=7$, then
$$
E_1=\{u_1u_{2},u_1u_{3},u_{1}u_{5},u_1u_{13},u_{1}u_{21},u_1u_{23},u_1u_{24}\}.
$$

In this case, the circulant graph $\Gamma_{24}$ can be found out;
see Figure 1 $(b)$.

\textbf{Step 2.} By the edge set $E_1$, we write the vector
$\alpha_{24}$ corresponding to $E_1$. If $|E_1|=9$, then
$$
\alpha_{24}=(0,1,1,0,1,0,1,0,0,0,0,0,1,0,0,0,0,0,1,0,1,0,1,1).
$$

If $|E_1|=7$, then
$$
\alpha_{24}=(0,1,1,0,1,0,0,0,0,0,0,0,1,0,0,0,0,0,0,0,1,0,1,1).
$$

\textbf{Step 3.} Change the first component of the vector
$\alpha_{24}$ to $\omega$. Denote by $\alpha'_{24}$ the new vector.
We generate a circulant matrix $A'_{24}$($A''_{24}$) from $\alpha'_{24}$($\alpha''_{24}$).
$$
\alpha_{24}'=(\omega,1,1,0,1,0,1,0,0,0,0,0,1,0,0,0,0,0,1,0,1,0,1,1)
$$
$$
\alpha_{24}''=(\omega,1,1,0,1,0,0,0,0,0,0,0,1,0,0,0,0,0,0,0,1,0,1,1).
$$

\textbf{Step 4.} By Algorithm $1$, we obtain the minimum distance
$d'_{24}=d''_{24}=8$. Therefore, both $\mathcal {C}_{24}'$ and
$\mathcal {C}_{24}''$ are proposed optimum codes over $GF(4)$.

The weight enumerators of the codes $\mathcal {C}_{24}'$ and
$\mathcal {C}_{24}''$ are
\begin{eqnarray*}
W_{\mathcal
{C}_{24}'}(z)&=&1+528z^{8}+13992z^{10}+171276z^{12}+1118040z^{14}+3773517z^{16}+6218520z^{18}\\
& &+4413948z^{20}+1034088z^{22}+33306z^{24}\\
\end{eqnarray*}
\begin{eqnarray*}
W_{\mathcal
{C}_{24}''}(z)&=&1+648z^{8}+13032z^{10}+174636z^{12}+1111320z^{14}+3781917z^{16}+6211800z^{18}\\
& &+4417308z^{20}+1033128z^{22}+33426z^{24}\\
\end{eqnarray*}

Applying the above method, we can obtain optimum graph codes over
$GF(4)$ from the first two generator vector.

\begin{tabbing}
\noindent\rule[0.25\baselineskip]{\textwidth}{1pt}\\

$n$~~~~~~$d$~~~~~~$L_n$~~~~~~First row of generator matrix\\

\noindent\rule[0.25\baselineskip]{\textwidth}{0.5pt}\\

$16$~~~~~~$6$~~~~~~$6$~~~~~~\textcolor{red}{$(\omega,1,1,0,1,0,0,0,1,0,0,0,1,0,1,1)$}\\

$22$~~~~~~$8$~~~~~~$8$~~~~~~\textcolor{red}{$(\omega,1,1,0,1,0,1,0,0,0,0,1,0,0,0,0,1,0,1,0,1,1)$}\\

$16$~~~~~~$6$~~~~~~$6$~~~~~~\textcolor{red}{$(\omega,1,0,1,0,1,0,0,1,0,0,1,0,1,0,1)$}\\

$19$~~~~~~$7$~~~~~~$7$~~~~~~\textcolor{red}{$(\omega,1,1,0,1,0,0,0,1,1,0,0,0,1,0,1,1)$}\\

\noindent\rule[0.25\baselineskip]{\textwidth}{1pt}
\end{tabbing}

\section{New codes from dense circulant graphs}

From the above section, we see that the circulant graphs under our
consideration are all relatively sparse. So one may think that only
sparse circulant graphs produce optimum graph codes. The following
fact gives it a negative answer.

Danielsen \cite{Danielsen} obtained an optimum additive code
$(30,2^{30},12)$ from the vector
$$
\alpha_{30}=(\omega,0,1,1,0,0,0,0,1,1,0,1,1,1,1,1,1,1,1,1,0,1,1,0,0,0,0,1,1,0),
$$
which corresponds to a circulant graph of order $30$ such that its
adjacent matrix $A_{30}$ is generated by the first row
$$
\alpha_{30}=(0,0,1,1,0,0,0,0,1,1,0,1,1,1,1,1,1,1,1,1,0,1,1,0,0,0,0,1,1,0).
$$
It is clear that the circulant graph is $17$-regular. Note that the
order of this graph is $30$, so the degree of each vertex is
relatively large, i.e., it is a relatively dense graph. Let

\[
\alpha_n=(\omega,0,1,1,0,0,0,0,1,1,0,\overbrace{1,1,\ldots,1,1}^{n-21},0,1,1,0,0,0,0,1,1,0).
\]

For $n=32$ and $n=34$, we have the following vectors,
$$
\alpha_{32}=(\omega,0,1,1,0,0,0,0,1,1,0,1,1,1,1,1,1,1,1,1,1,1,0,1,1,0,0,0,0,1,1,0),
$$
$$
\alpha_{34}=(\omega,0,1,1,0,0,0,0,1,1,0,1,1,1,1,1,1,1,1,1,1,1,1,1,0,1,1,0,0,0,0,1,1,0).
$$
By Algorithm $1$, both $\mathcal {C}_{32}$ and $\mathcal {C}_{34}$
are proposed optimum codes over $GF(4)$.

The weight enumerators of the codes $\mathcal {C}_{32}$ and
$\mathcal {C}_{34}$ are
\begin{eqnarray*}
W_{\mathcal
{C}_{32}}(z)&=&1+1325z^{10}+41973z^{12}+745155z^{14}+8030541z^{16}+53150370z^{18}+213875634z^{20}\\
& &+510617670z^{22}+691665390z^{24}+491629473z^{26}+159600905z^{28}+17838471z^{30}\\
& &+286740z^{32}\\
\end{eqnarray*}
\begin{eqnarray*}
W_{\mathcal
{C}_{34}}(z)&=&1+492z^{10}+14373z^{12}+291849z^{14}+3494061z^{16}+26279603z^{18}+123536402z^{20}\\
& &+357928154z^{22}+620714798z^{24}+614055698 z^{26}+319190777z^{28}+75747789z^{30}\\
& &+6157556z^{32}+72095z^{34}\\
\end{eqnarray*}


\begin{thebibliography}{11}

\bibitem{Bermond}
J. C. Bermond, F. Comellas, D. F. Hsu,
\emph{Distributed loop computer networks: A survey}, J. Parallel
Distributed Comput. 24 (1995), 2-10.

\bibitem{Boesch}
F. T. Boesch, J. F. Wang, \emph{Reliable circulant networks with
minimum transmission delay}, IEEE Trans. Circuits Syst. 32(1985),
1286-1291.

\bibitem{Bondy}
J. A. Bondy, U. S. R. Murty,
{\it Graph Theory}, GTM 244, Springer, 2008.


\bibitem{Danielsen}
L. E. Danielsen, \emph{On Self-Dual Quantum Codes},
Graphs and Boolean Functions, 2005.

\bibitem{Danielsen2}
L. E. Danielsen, \emph{Graph-based classification of self-dual
additive codes over finite field}, Adv. Math. Commun. 3(4)(2009),
329-348.

\bibitem{Danielsen3}
L. E. Danielsen, \emph{On the Classification of Hermitian Self-Dual Additive Codes over
GF(9)}, IEEE Trans. Inform. Theory 58(8)(2012), 5500-5511.

\bibitem{Danielsen4}
L. E. Danielsen, M. G. Parker, \emph{Directed Graph
Representation of Half-Rate Additive Codes over GF(4)}, Des. Codes
Cryptogr. 59(2011), 119-130.

\bibitem{DP}
L. E. Danielsen, M. G. Parker,
\emph{On the classification of all self-dual additive codes over
$GF(4)$ of length up to $12$}, J. Combin. Theory, Series $A$,
\textbf{113}(2006), 1351-1367.

\bibitem{GGMG}
D. G. Glynn, T. A. Gulliver, J. G. Marks, M. K. Gupta,
\emph{The Geometry of Additive Quantum Codes}, Preface, Springer,
2006.

\bibitem{HP}
W. C. Huffman, Vera Pless,
{\it Fundamentals of Error-Correcting Codes}, Cambridge University,
2003.

\bibitem{Mans}
B. Mans, F. Pappalardi, I. Shparlinski, \emph{On the spectral Adam
property for circulant graphs}, Discrete Math. 254(1-3)(2002),
309-329.

\bibitem{Meijer}
P. T. Meijer,
{\it Connectivities and Diameters of Circulant Graps},
B.Sc.(Honors), Simon Fraser University, 1987.

\bibitem{Monakhova}
E. A. Monakhova, \emph{A survey on undirected
circulant graphs}, Discrete Mathematics, Algorithms and
Applications, 4(1)(2012), DOI: 10.1142/S1793830912500024.


\bibitem{Muzychuk}
M. E. Muzychuk, G. Tinhofer, \emph{Recognizing circulant graphs of
prime order in polynomial time}, Electron. J. Combin. 5(1)(1998),
501-528.

\bibitem{Stein}
W. A. Stein et al., \emph{Sage Mathematics Software (Version
  6.1.1)}, The Sage Development Team, 2014, http://www.sagemath.org.

\bibitem{Tonchev}
V. Tonchev,
\emph{Error-correcting codes from graphs}, Disrete Math.,
\textbf{257}(2002), 549-557.

\bibitem{Varbanov}
Z. Varbanov,
\emph{Additive circulent graph codes over GF(4)}, Math. Maced.
6(2008), 73-79.



\end{thebibliography}
\end{document}